\documentclass[10pt]{article}
\usepackage{epsf}
\usepackage{amsmath}

\allowdisplaybreaks

\usepackage[showframe=false]{geometry}
\usepackage{changepage}

\usepackage{epsfig}
\usepackage{amssymb}

\usepackage{amsthm}
\usepackage{setspace}
\usepackage{cite}

\usepackage{mcite}

\usepackage{algorithmic}  
\usepackage{algorithm}

\usepackage{shadow}
\usepackage{fancybox}
\usepackage{fancyhdr}

\usepackage{color}
\usepackage[usenames,dvipsnames,svgnames,table]{xcolor}
\newcommand{\bl}[1]{\textcolor{blue}{#1}}

\definecolor{mypurple}{rgb}{.4,.0,.5}

\usepackage[hyphens]{url}

\usepackage[colorlinks=true,
            linkcolor=black,
            urlcolor=blue,
            citecolor=purple]{hyperref}

\usepackage{breakurl}
\def\w{{\bf w}}

\def\y{{\bf y}}

\def\x{{\bf x}}

\def\x{{\mathbf x}}

\def\w{{\bf w}}

\def\x{{\bf x}}
\def\y{{\bf y}}

\def\f{{\bf f}}

\def\diag{{\rm diag}\,}

\def\be{\begin{equation}}
\def\ee{\end{equation}}
\def\ba{\left[\begin{array}}
\def\ea{\end{array}\right]}

\def\w{{\bf w}}

\def\x{{\bf x}}
\def\y{{\bf y}}

\def\f{{\bf f}}

\def\1{{\bf 1}}
\def\oness{{\bf 1}^{(s)}}

\def\g{{\bf g}}
\def\0{{\bf 0}}

\def\erfc{\mbox{erfc}}
\def\erfi{\mbox{erfi}}
\def\diag{\mbox{diag}}
\def\keta{k_{\eta}}

\def\mR{{\mathbb R}}
\def\mN{{\mathbb N}}
\def\mE{{\mathbb E}}

\def\phiint{\phi_{int}}
\def\phiext{\phi_{ext}}

\def\calF{{\cal F}}

\def\lp{\left (}
\def\rp{\right )}

\newtheorem{theorem}{Theorem}

\begin{document}

\begin{singlespace}

\title {Partial $\ell_1$ optimization in random linear systems -- finite dimensions 
}
\author{
\textsc{Mihailo Stojnic
\footnote{e-mail: {\tt flatoyer@gmail.com}} }}
\date{}
\maketitle

\centerline{{\bf Abstract}} \vspace*{0.1in}

In this paper we provide a complementary set of results to those we present in our companion work \cite{Stojnicl1HidParasymldp} regarding the  behavior of the so-called partial $\ell_1$ (a variant of the standard $\ell_1$ heuristic often employed for solving under-determined systems of linear equations). As is well known through our earlier works \cite{StojnicICASSP10knownsupp,StojnicTowBettCompSens13}, the partial $\ell_1$ also exhibits the phase-transition (PT) phenomenon, discovered and well understood in the context of the standard $\ell_1$ through Donoho's and our own works \cite{DonohoPol,DonohoUnsigned,StojnicCSetam09,StojnicUpper10}. \cite{Stojnicl1HidParasymldp} goes much further though and, in addition to the determination of the partial $\ell_1$'s phase-transition curves (PT curves) (which had already been done in \cite{StojnicICASSP10knownsupp,StojnicTowBettCompSens13}), provides a substantially deeper understanding of the PT phenomena through a study of the underlying large deviations principles (LDPs). As the PT and LDP phenomena are by their definitions related to large dimensional settings, both sets of our works, \cite{StojnicICASSP10knownsupp,StojnicTowBettCompSens13} and \cite{Stojnicl1HidParasymldp}, consider what is typically called the asymptotic regime. In this paper we move things in a different direction and consider finite dimensional scenarios. Basically, we provide explicit performance characterizations for any given collection of systems/parameters dimensions. We do so for two different variants of the partial $\ell_1$, one that we call exactly the partial $\ell_1$ and another one, possibly a bit more practical, that we call the hidden partial $\ell_1$. Finally, we also show for both of these variants how one can bridge between the finite dimensional settings considered here and the infinite dimensional ones considered in  \cite{Stojnicl1HidParasymldp} (and earlier in \cite{StojnicICASSP10knownsupp,StojnicTowBettCompSens13} as well).

\vspace*{0.25in} \noindent {\bf Index Terms: Finite dimensions;
partial $\ell_1$; linear systems of equations; sparse solutions}.

\end{singlespace}

\section{Introduction}
\label{sec:back}

In this paper we will study a well known topic of linear systems with sparse solutions. Obviously this is a very large area and we focus on several particular problems that relate to the so-called $\ell_1$ heuristic often used for their solving. As the descriptions of the problems that we will study here will closely resemble the description of the problems studied in \cite{Stojnicl1HidParasymldp}, we will refrain from repeating portions of the discussion from \cite{Stojnicl1HidParasymldp} that relate to the importance of these problems and relevant prior work (instead, we will assume a solid level of familiarity with these and, in general, with what is done in \cite{Stojnicl1HidParasymldp}).

To introduce the problems that will be of interest here, we start with the introduction of the standard $\ell_1$. Let $A$ be an $m\times n$ ($m\leq n$) dimensional real matrix (which we will often refer to as the system matrix) and let $\tilde{\x}$ be an $n$ dimensional real vector (for short we say, $A\in \mR^{m\times n}$ and $\tilde{\x}\in \mR^{n}$). $\tilde{\x}$ will be called $k$-sparse if no more than $k$ of its entries are different from zero. Let $\y$ be such that
\begin{equation}
\y\triangleq A\tilde{\x}. \label{eq:defy}
\end{equation}
Assuming that $A$ and $\y$ in (\ref{eq:defy}) are given, one would like to deteremine/recover $\tilde{\x}$. Or in other words, can one find the unknown $\x$ in the following system of linear equations
\begin{equation}
A\x=\y. \label{eq:system}
\end{equation}
An under-determined regime ($m\leq n$) will be of our interest here. Moreover, $\tilde{\x}$ will be assumed to be $k$ sparse, and the algebraic properties of $A$ and the relation between $k$, $m$, and $n$ will be assumed such that the $k$-sparse $\x$ that solves (\ref{eq:system}) is unique and that there is no sparser $\x$ that solves (\ref{eq:system}) (either deterministically or when $A$ is random statistically). Then one often rewrites (\ref{eq:system}) as
\begin{eqnarray}
\mbox{min} & & \|\x\|_{0}\nonumber \\
\mbox{subject to} & & A\x=\y, \label{eq:l0}
\end{eqnarray}
where $\|\x\|_{0}$ essentially counts the number of the nonzero entries of $\x$. Finding the sparsest $\x$ in (\ref{eq:l0}) (which we will technically call solving (\ref{eq:l0})) is not an easy task (see, e.g. \cite{StojnicCSetam09,Stojnicl1RegPosasymldp,StojnicReDirChall13}). Many heuristics have been developed over last several decades though (see, e.g. \cite{JATGomp,NeVe07,DTDSomp,NT08,DaiMil08,DonMalMon09}) that often achieve a solid level of success in solving (\ref{eq:l0}). As mathematically the strongest, we view, the following $\ell_1$-optimization relaxation of (\ref{eq:l0})
\begin{eqnarray}
\mbox{min} & & \|\x\|_{1}\nonumber \\
\mbox{subject to} & & A\x=\y. \label{eq:l1}
\end{eqnarray}

\subsection{Partial $\ell_1$ -- definitions}
\label{sec:parl1}

Various forms of performance characterization (PTs, LDPs, and even finite dimension success/failure rates) of the above $\ell_1$ in statistical scenarios have been fully settled through \cite{DonohoPol,DonohoUnsigned,StojnicCSetam09,StojnicUpper10,Stojnicl1RegPosasymldp,Stojnicl1RegPosfinn}. All of these works typically produce the following two major reasons for a large popularity of the above $\ell_1$ heuristic: 1) it is a simple linear program and 2) it has provably excellent performance abilities when used for recovering the unknown $\x$ in
(\ref{eq:l0}). Naturally the following question then arises: can one substantially improve on either of these points without substantially changing the other? In other words, can one design an algorithm that is at least as fast as the standard $\ell_1$ and has performance abilities provably substantially better than those established in \cite{DonohoPol,DonohoUnsigned,StojnicCSetam09,StojnicUpper10,Stojnicl1RegPosasymldp,Stojnicl1RegPosfinn} or can one design an algorithm that is substantially faster than the standard $\ell_1$ and has performance abilities provably at least as good as those established in \cite{DonohoPol,DonohoUnsigned,StojnicCSetam09,StojnicUpper10,Stojnicl1RegPosasymldp,Stojnicl1RegPosfinn}. An interesting line of work that can be utilized to initiate thinking about these improvements, was started in \cite{VasLu09}. Namely, in \cite{VasLu09} and later on in our works \cite{Stojnicl1HidParasymldp,StojnicICASSP10knownsupp,StojnicTowBettCompSens13} (as well as in \cite{SteChr09} to a degree) the following modification of the standard $\ell_1$ from (\ref{eq:l1}) was considered
\begin{eqnarray}
\mbox{min} & & \sum_{i\notin\Pi} |\x_i|\nonumber \\
\mbox{subject to} & & A\x=\y, \label{eq:l1imp}
\end{eqnarray}
where $\Pi$ is a set of cardinality $k-k_{\eta}$ (clearly, $k_{\eta}< k$). This modification, which we will often refer to as the \emph{partial} $\ell_1$, is originally developed for the scenarios where the portion of the support of the unknown $\x$ in (\ref{eq:system}) (which is the set of the nonzero locations of $\x$ and which we denote by $supp(\x)$) is beforehand known. This portion is, of course, exactly the set $\Pi$. However, this modification is also at the heart of the so-called iterative $\ell_1$ strategies (see, e.g. \cite{SteChr09}) that attempt to improve on the standard $\ell_1$. In \cite{StojnicICASSP10knownsupp,StojnicTowBettCompSens13}, we managed to fully characterize the partial $\ell_1$'s phase-transitions which essentially explained what kind of gain one can expect asymptotically from knowing $\Pi$, i.e. from knowing a portion of the support of the unknown $\x$ in (\ref{eq:system}), basically $supp(\tilde{\x})$. A much deeper understanding of the partial $\ell_1$ behavior in the entire transition zone is provided in \cite{Stojnicl1HidParasymldp}. Not only do \cite{StojnicICASSP10knownsupp,StojnicTowBettCompSens13,Stojnicl1HidParasymldp} establish the asymptotically precise gains that (\ref{eq:l1imp}) makes over (\ref{eq:l1}) in scenarios where a feedback in the form of a known $\Pi$ is available, they also suggest that the original $\ell_0$ problem from (\ref{eq:l0}) can be transformed into a partial recovery problem where one searches for sets $\Pi$ instead of searching for the entire $supp(\tilde{\x})$. Of course things are a bit more involved as one has to be careful as to what kind of sets $\Pi$ he deals with. This discussion goes of course beyond what our main interests of study here are and we instead of pursuing a further consideration in this direction refer to \cite{StojnicICASSP10knownsupp,StojnicTowBettCompSens13,Stojnicl1HidParasymldp} where more on this can be found. Here though, we focus on the resulting mathematical problems that were the subject of interest in \cite{StojnicICASSP10knownsupp,StojnicTowBettCompSens13,Stojnicl1HidParasymldp}. Differently from \cite{StojnicICASSP10knownsupp,StojnicTowBettCompSens13,Stojnicl1HidParasymldp} where the nature of the PT and LDP phenomena basically implied an asymptotic performance analysis of (\ref{eq:l1imp}), we here provide its a finite dimensional counterpart. The results that we will present below will be in flavor similar to those that we presented in the introductory paper \cite{Stojnicl1RegPosfinn}. In fact, to be a bit more precise, they will relate to the asymptotic results of \cite{StojnicICASSP10knownsupp,StojnicTowBettCompSens13,Stojnicl1HidParasymldp}, in pretty much the same way the finite dimensional considerations from \cite{Stojnicl1RegPosfinn} relate to the asymptotic ones from \cite{StojnicCSetam09,StojnicUpper10,Stojnicl1RegPosasymldp}.

We will organize the paper so that the presentation is split into two main sections. In the first one we will consider the partial $\ell_1$ from (\ref{eq:l1imp}). In the second one we will discuss its a possibly more practical variant that we will refer to as the \emph{hidden partial} $\ell_1$. Before switching to a detailed discussion about the partial $\ell_1$ we will just briefly introduce the hidden partial $\ell_1$.

\subsection{Hidden partial $\ell_1$ -- definitions}
\label{sec:hidl1}

As discussed above, the partial $\ell_1$ can be utilized not only when there is an available feedback about the (portion of) $supp(\tilde{\x})$, but also as a part of a general strategy to improve on the standard $\ell_1$. Our results \cite{StojnicICASSP10knownsupp,StojnicTowBettCompSens13} in fact provided a fully mathematically rigorous, asymptotically exact gain that the partial $\ell_1$ achieves over the standard $\ell_1$. As discussed in \cite{Stojnicl1HidParasymldp}, the results that relate to the partial $\ell_1$ of course assume that an a priori available set $\Pi$ contains nothing more than a subset of $supp(\x)$. To handle scenarios where $\Pi$ contains so to say a bit of ``noise" (i.e. some of the locations of zeros of $\x$), we in \cite{StojnicTowBettCompSens13} introduced the hidden partial $\ell_1$ modification of the partial $\ell_1$. To fully explain this modification we will need a couple of additional definitions introduced in \cite{StojnicTowBettCompSens13}. Let $\kappa\subset \{1,2,\dots,n\}$ and let the cardinality of $\kappa$ be $k$ (we will for the simplicity choose $k$; however our results easily extend to any other value). Let $\Pi$ be the intersection of the set of nonzero locations of $\x$ ($supp(\x)$) and $\kappa$. As earlier, $\Pi$ is the set that is known to contain locations of some of the nonzero elements of $\x$. Differently though from what was the case earlier, $\Pi$ is not known now. Instead $\kappa$ is now known and the fact that $\Pi\in\kappa$. For the concreteness, we will again assume that the cardinality of $\Pi$ is $k_{\eta}$ (where $k_{\eta}<k$) and will view $\x$ as a vector with \emph{hidden} partially known support (clearly, $\kappa$ will represent the estimate of $\x$'s support ($supp(\x)$)). Then the above mentioned hidden partial $\ell_1$ assumes the following slight adjustment to (\ref{eq:l1imp})
\begin{eqnarray}
\mbox{min} & & \sum_{i\notin\kappa} |\x_i|\nonumber \\
\mbox{subject to} & & A\x=\y. \label{eq:l1imphidden}
\end{eqnarray}
It is rather obvious that this algorithm will be of use in scenarios where there is an available decent estimate $\kappa$ of $supp(\x)$. Why would this be a more practical variant of the partial $\ell_1$? If one views the partial $\ell_1$ as a tool that can be used to improve on the standard $\ell_1$ one then needs to ``feed" (\ref{eq:l1imp}) with $\Pi$. However, it is much easier to determine a $\kappa$ than a $\Pi$. In fact, many heuristics that attempt to solve (\ref{eq:system}) actually do so. Namely, even when they fail to correctly uncover the entire $supp(\tilde{\x})$, they still output an estimate of it that more often than not contains a large number of the elements of $supp(\tilde{\x})$ (that is precisely what $\kappa$ can be thought of as well).


\section{Partial $\ell_1$ -- finite dimensional analysis}
\label{sec:posl1}

In this section we analyze the partial $\ell_1$ from (\ref{eq:l1imp}). We will try to follow as much as possible the corresponding analysis for the standard $\ell_1$ presented in \cite{Stojnicl1RegPosfinn}. To that end, we will also try to avoid repeating all the arguments that remain the same and instead will focus on emphasizing the ones that are different. As usual, (and as was done in \cite{StojnicCSetam09,StojnicUpper10,Stojnicl1RegPosasymldp,Stojnicl1HidParasymldp,Stojnicl1RegPosfinn}), we start things off by recalling on a couple of results that we established in \cite{StojnicICASSP10knownsupp,StojnicTowBettCompSens13}, (these were, of course, the key components of the analysis done in \cite{StojnicICASSP10knownsupp,StojnicTowBettCompSens13}; their standard $\ell_1$ counterparts were as important in \cite{StojnicCSetam09,StojnicUpper10,Stojnicl1RegPosasymldp,Stojnicl1RegPosfinn} as well as in a large sequence of our work that followed later on).

For the concreteness of the exposition and without loss of generality we will assume that the elements $\x_{k+1},\x_{k+2},\dots,\x_{n}$ of $\x$ are equal to zero and that the elements $\x_{1},\x_{2},\dots,\x_k$ have fixed signs, say all positive (this is of course not known beforehand and cannot be utilized in the design of algorithms). Also, for concreteness and without loss of generality, let $\Pi=\{1,2,\dots,\keta\}$. The following is then a partial $\ell_1$ adaptation of the result proven for the general $\ell_1$ in \cite{StojnicCSetam09,StojnicUpper10,StojnicICASSP09} and, as mentioned above, is among the key unsung heros that enabled running the entire machinery developed overthere.
\begin{theorem}(\cite{StojnicICASSP10knownsupp,StojnicTowBettCompSens13} Nonzero elements of $\x$ have fixed signs and location)
Assume that an $m\times n$ system matrix $A$ is given. Let $\x$
be a $k$ sparse vector. Also let $\x_{k+1}=\x_{k+2}=\dots=\x_{n}=0$. Let the signs of $\x_{1},\x_{2},\dots,\x_k$ be fixed, say all positive. Also, let $\Pi=\{1,2,\dots,\keta\}$. Further, assume that $\y=A\x$ and that $\w$ is
a $n\times 1$ vector. If
\begin{equation}
(\forall \w\in \textbf{R}^n | A\w=0) \quad  -\sum_{i=\keta+1}^{k} \w_i<\sum_{i=k+1}^{n}|\w_{i}|,
\label{eq:posthmcond1}
\end{equation}
then the solutions of (\ref{eq:l0}) and (\ref{eq:l1imp}) coincide. Moreover, if
\begin{equation}
(\exists \w\in \textbf{R}^n | A\w=0) \quad  -\sum_{i=\keta+1}^{k} \w_i\geq \sum_{i=k+1}^{n}|\w_{i}|,
\label{eq:posthmcond2}
\end{equation}
then there is an $\x$ from the above set of $\x$'s with fixed location of nonzero elements such that the solution of (\ref{eq:l0}) is not the solution of (\ref{eq:l1imp}).
\label{thm:posthmregposcond}
\end{theorem}
To facilitate the exposition we set
\begin{equation}
C^{(p)}_w\triangleq C^{(p)}_w(k,m,n,\keta)=\{\w\in \mR^n| \quad -\sum_{i=\keta+1}^k \w_i\geq \sum_{i=k+1}^{n}|\w_{i}|\}.\label{eq:defSw}
\end{equation}
$C^{(p)}_w$ is a polyhedral cone (from this point on we will often assume a decent level of familiarity with some of the well known concepts in high-dimensional geometry; more on them though can be found in e.g. \cite{BG,SantaloBookIntGeom76}). Following what was done in \cite{Stojnicl1RegPosfinn} we can write
\begin{equation}
p^{(p)}_{err}(k,m,n,k_{\eta})=P(G^{(sub)}\cap C^{(p)}_w\neq \emptyset)=2\sum_{l=m+2j+1,j\in \mN_0}^{n} \sum_{F^{(l)}\in \calF^{(l)}}\phiint(0,F^{(l)})\phiext(F^{(l)},C^{(p)}_w),\label{eq:anal3}
\end{equation}
where $p^{(p)}_{err}(k,m,n,k_{\eta})$ is the probability of error that the $k$-sparse solution of (\ref{eq:l0}) is not the solution of (\ref{eq:l1imp}), $F^{(l)}$ are the $l$-faces of $C^{(p)}_w$, $\calF^{(l)}$ is the set of all $l$-faces of $C^{(p)}_w$, and $\phiint(\cdot,\cdot)$ and $\phiext(\cdot,\cdot)$ are the so-called internal and external angles, respectively (see, e.g. \cite{Santalo,PMM,AS,BG}). Now, (\ref{eq:anal3}) is a nice conceptual solution to the problem of determining $p^{(p)}_{err}(k,m,n,k_{\eta})$. To have (\ref{eq:anal3}) be fully operational though one would need to be able to handle the angles $\phiint(\cdot,\cdot)$ and $\phiext(\cdot,\cdot)$ which is typically a very hard problem and very rarely solvable. In \cite{Stojnicl1RegPosfinn} we presented a fairly elegant way to compute these angles and below we discuss how it can be utilized for the problems of interest here. We will first deal separately with the internal angles and then afterwards switch to the external angles.

\subsection{Internal angles}
\label{sec:intang}

In this section we analyze the internal angles appearing in (\ref{eq:anal3}), i.e. $\phiint(0,F^{(l)})$. Similarly to what was done in \cite{Stojnicl1RegPosfinn}, we will distinguish between two cases: 1) $l<n$ and 2) $l=n$. For $l<n$, we have for the set of all $l$-faces $\calF^{(l,p)}_1$
\begin{multline}
\calF^{(l,p)}_1 \triangleq\{\w\in \mR^n| \quad -\sum_{i=\keta+1}^k \w_i= (\oness)^T\w_{I_r\setminus I^{(l)}_1},\diag(\oness)\w_{I_r\setminus I^{(l)}_1}\geq 0,\w_{I^{(l)}_1}=0,\\I^{(l)}_1\subset I_r,|I^{(l)}_1|=n-l-1,\oness\in\{-1,1\}^{l-k+1}\}.\label{eq:intanal2}
\end{multline}
The cardinality of set $\calF^{(l,p)}_1$ is easily given by
\begin{equation}
c^{(l)}_1\triangleq |\calF^{(l,p)}_1|=2^{l-k+1}\binom{n-k}{n-l-1}, l\in \{k-1,k,\dots,n-1\}.\label{eq:intanal4}
\end{equation}
(\ref{eq:anal3}) can then be rewritten in the following way
\begin{eqnarray}
p^{(p)}_{err}(k,m,n,k_{\eta}) & = & 2\sum_{l=m+2j+1,j\in \mN_0}^{n}  \sum_{F^{(l,p)}_1\in \calF^{(l,p)}_1}\phiint(0,F^{(l,p)}_1)\phiext(F^{(l,p)}_1,C^{(p)}_w)\nonumber \\
& = & 2 ( \sum_{l=m+2j+1,j\in \mN_0}^{n-1} c^{(l)}_1\phiint(0,F^{(l,p)}_1)\phiext(F^{(l,p)}_1,C^{(p)}_w) +  \phiint(0,C^{(p)}_w)\phiext(C^{(p)}_w,C^{(p)}_w)),\nonumber \\
\label{eq:intanal6}
\end{eqnarray}
where due to symmetry $F^{(l,p)}_1$ is basically any of the elements from set $\calF^{(l,p)}_1$. For the concreteness we choose $I^{(l)}_1=\{l+2,l+3,\dots,n\}$ and $\oness=\1_{(l-k+1)\times 1}$ and consequently have
\begin{equation}
F^{(l,p)}_1\triangleq F^{(l,p)}_1(k,\keta) =\{\w\in \mR^n| \quad -\sum_{i=\keta+1}^k \w_i= \sum_{i=k+1}^{l+1}\w_{i},\w_{k+1:l+1}\geq 0,\w_{l+2:n}=0\}, l\in \{k-1,k,\dots,n-1\}.\label{eq:intanal7}
\end{equation}
Below we separately discuss $\phiint(0,F^{(l,p)}_1)$ and $\phiint(0,C^{(p)}_w)$.

\subsubsection{Determining $\phiint(0,F^{(l,p)}_1)$}
\label{sec:int1ang}

To compute $\phiint(0,F^{(l,p)}_1)$ we will follow the ``Gaussian coordinates in an orthonormal basis" strategy presented in \cite{Stojnicl1RegPosfinn}. We recall that the strategy assumes two steps: 1) Finding an orthonormal basis in the subspace where the angle is being computed; 2) Expressing the content of the angle in terms of the Gaussian coordinates of the computed orthonormal basis.

1) Similarly to what was done in \cite{Stojnicl1RegPosfinn}, for the orthonormal basis we will use the column vectors of the following matrix
\begin{equation}
B_{int,1}=\begin{bmatrix}
            \begin{bmatrix}
              B \\
              \0_{1\times (l-1)}
            \end{bmatrix}  & \begin{bmatrix}
            \0_{\keta\times 1}\\
              -\1_{(l-\keta)\times 1} \\
              l-\keta
            \end{bmatrix}\frac{1}{\sqrt{(l-\keta)^2+l-\keta}}\\
            \0_{(n-l-1)\times (l-1)} & \0_{(n-l-1)\times 1}
          \end{bmatrix},\label{eq:posint1anal1}
\end{equation}
where $B$ is an $l\times (l-1)$ orthonormal matrix such that $[\begin{bmatrix}
\0_{1\times \keta} & \1_{1\times (l-\keta)}\end{bmatrix} B=\0_{(l-1)\times 1}$, and $\1$ and $\0$ are matrices of all ones or zeros, respectively of the sizes given in their indexes. One can easily confirm that $B_{int,1}^TB_{int,1}=I$. Moreover, $F^{(l,p)}_1$ is indeed in the subspace spanned by the columns of $B_{int,1}$ since its normal vector
$\f^{(l,1)}=\begin{bmatrix}\0_{1\times \keta} & -\1_{1\times (l+1-\keta)} & \0_{1\times (n-l-1)}\end{bmatrix}^T$ does satisfy $(\f^{(l,1)})^T B_{int,1}=\0_{1\times l}$.

2) Following further \cite{Stojnicl1RegPosfinn}, we will operate in this orthonormal basis through the standard normal (i.e. Gaussian) coordinates. In other words, $\g\in\mR^{l}$ will be assumed to have $l$ i.i.d standard normal components and only those $\g$ for which $B_{int,1}\g\in F^{(l,p)}_1$ will be allowed. This means that we will be interested in the following set of $\g$'s, $G^{(l)}_1$,
\begin{equation}
G^{(l)}_1 =\{\g\in \mR^l| \w=B_{int,1}\g \quad \mbox{and}\quad \w_i\geq 0,k+1\leq i\leq l+1\}, l\in \{k-1,k,\dots,n-1\}.\label{eq:posint1anal2}
\end{equation}
Relying on (\ref{eq:posint1anal1}), one can reformulate (\ref{eq:posint1anal2}) for any $l\in \{k-1,k,\dots,n-1\}$ in the following way
\begin{equation}
G^{(l)}_1 =\{\g\in \mR^l| \w_{1:(l+1)}=\begin{bmatrix}
            \begin{bmatrix}
              B \\
              \0_{1\times (l-1)}
            \end{bmatrix}  & \begin{bmatrix}
              \0_{\keta\times 1} \\
              -\1_{(l-\keta)\times 1} \\
              l-\keta
            \end{bmatrix}\frac{1}{\sqrt{(l-\keta)^2+l-\keta}}
          \end{bmatrix}\g \quad \mbox{and}\quad \w_i\geq 0,k+1\leq i\leq l+1\}.\label{eq:posint1anal3}
\end{equation}
By the definition of the internal angle we then have
\begin{equation}
\phiint(0,F^{(l,p)}_1)=\frac{1}{(2\pi)^{\frac{l}{2}}}\int_{G^{(l)}_1}e^{-\frac{\g^T\g}{2}}d\g.\label{eq:posint1anal4}
\end{equation}
The following change of variables (clearly relying on (\ref{eq:posint1anal3}))
\begin{equation}
\w_{1:l}=\begin{bmatrix}
              B  &  \begin{bmatrix}
                      \0_{\keta\times 1} \\
                      -\1_{(l-\keta)\times 1}
                    \end{bmatrix}\frac{1}{\sqrt{(l-\keta)^2+l-\keta}}
          \end{bmatrix}\g.
          \label{eq:posint1anal5}
\end{equation}
easily gives
\begin{equation}
\g=\begin{bmatrix}
              B^T  \\  \begin{bmatrix} \0_{1\times \keta} & -\1_{1\times (l-\keta)}\end{bmatrix}\sqrt{\frac{(l-\keta)^2+l-\keta}{(l-\keta)^2}}
          \end{bmatrix}\w_{1:l},
          \label{eq:posint1anal6}
\end{equation}
and ultimately
\begin{eqnarray}
\g^T\g & = & \w_{1:l}^T\begin{bmatrix}
              B^T  \\  \begin{bmatrix} \0_{1\times \keta} & -\1_{1\times (l-\keta)}\end{bmatrix}\sqrt{\frac{(l-\keta)^2+l-\keta}{(l-\keta)^2}}
          \end{bmatrix}^T\begin{bmatrix}
              B^T  \\  \begin{bmatrix} \0_{1\times \keta} & -\1_{1\times (l-\keta)}\end{bmatrix}\sqrt{\frac{(l-\keta)^2+l-\keta}{(l-\keta)^2}}
          \end{bmatrix}\w_{1:l}\nonumber \\
          &= &\w_{1:l}^T\begin{bmatrix}
              B^T  \\  \begin{bmatrix} \0_{1\times \keta} & -\1_{1\times (l-\keta)}\end{bmatrix}\sqrt{\frac{l-\keta}{(l-\keta)^2}}
          \end{bmatrix}^T\begin{bmatrix}
              B^T  \\  \begin{bmatrix} \0_{1\times \keta} & -\1_{1\times (l-\keta)}\end{bmatrix}\sqrt{\frac{l-\keta}{(l-\keta)^2}}
          \end{bmatrix}\w_{1:l}\nonumber \\
          & &+\w_{1:l}^T\begin{bmatrix} \0_{1\times \keta} & -\1_{1\times (l-\keta)}\end{bmatrix}^T
             \begin{bmatrix} \0_{1\times \keta} & -\1_{1\times (l-\keta)}\end{bmatrix}\w_{1:l}\nonumber \\
             & = & \w_{1:l}^T \w_{1:l}+\w_{1:l}^T\begin{bmatrix} \0_{1\times \keta} & -\1_{1\times (l-\keta)}\end{bmatrix}^T
             \begin{bmatrix} \0_{1\times \keta} & -\1_{1\times (l-\keta)}\end{bmatrix}\w_{1:l}.
          \label{eq:posint1anal7}
\end{eqnarray}
A combination of(\ref{eq:posint1anal3}) and (\ref{eq:posint1anal5}) also easily implies that $\w_{k+1:l}\geq 0$. Moreover, from (\ref{eq:posint1anal3}) we have $\g_l\geq 0$ which, by (\ref{eq:posint1anal6}), means that $\begin{bmatrix} \0_{1\times \keta} & -\1_{1\times (l-\keta)}\end{bmatrix}\w_{1:l}\geq 0$. It is also not that hard to compute the following Jacobian of the above change of variables in (\ref{eq:posint1anal5})
\begin{equation}
J=\frac{1}{\sqrt{\begin{bmatrix}
              B  &  \begin{bmatrix}
                      \0_{\keta\times 1} \\
                      -\1_{(l-\keta)\times 1}
                    \end{bmatrix}\frac{1}{\sqrt{(l-\keta)^2+l-\keta}}
          \end{bmatrix}^T\begin{bmatrix}
              B  &  \begin{bmatrix}
                      \0_{\keta\times 1} \\
                      -\1_{(l-\keta)\times 1}
                    \end{bmatrix}\frac{1}{\sqrt{(l-\keta)^2+l-\keta}}
          \end{bmatrix}}}=\sqrt{l+1-\keta}.\label{eq:posint1anal7a}
\end{equation}
Instead of (\ref{eq:posint1anal4}) we can now write
\begin{multline}
\phiint(0,F^{(l,p)}_1)   =  \frac{J}{(2\pi)^{\frac{l}{2}}}\int_{-\1_{1\times (l-\keta)}\w_{\keta+1:l}\geq 0,\w_{k+1:l}\geq 0}e^{-\frac{\w_{1:l}^T\w_{1:l}+\w_{\keta+1:l}^T( -\1_{1\times (l-\keta)})^T
             (-\1_{1\times (l-\keta)})\w_{\keta+1:l}}{2}}d\w_{1:l}\\
 = \frac{J}{(2\pi)^{\frac{l-\keta}{2}}}\int_{-\1_{1\times (l-\keta)}\w_{\keta+1:l}\geq 0,\w_{k+1:l}\geq 0}e^{-\frac{\w_{\keta+1:l}^T\w_{\keta+1:l}+\w_{\keta1:l}^T( -\1_{1\times (l-\keta)})^T
             (-\1_{1\times (l-\keta)})\w_{\keta1:l}}{2}}d\w_{\keta+1:l}.\label{eq:posint1anal8}
\end{multline}
We then recognize that the above form is exactly the same as the one obtained for the corresponding internal angle in the case of the standard $\ell_1$ studied in \cite{Stojnicl1RegPosfinn}. To be completely precise the following change $l\leftarrow l-\keta$ and $k\leftarrow k-\keta$ makes the angles identical. One can then skip the calculation of $\phiint(0,F^{(l,p)}_1)$ and instead use already computed corresponding quantity in \cite{Stojnicl1RegPosfinn} with the adjusted values for $l$ and $k$. One finally obtains
\begin{eqnarray}
\phiint(0,F^{(l,p)}_1)=   \frac{\sqrt{l+1-\keta}}{2^{l-k+1}\sqrt{2\pi}}\int_{-\infty}^{\infty} \lp1-i\erfi\lp\frac{t}{\sqrt{2}}\rp\rp^{l-k+1} e^{-\frac{(l+1-\keta)t^2}{2}}dt.
\label{eq:int1anal11}
\end{eqnarray}

\subsubsection{Determining $\phiint(0,C^{(p)}_w)$}
\label{sec:int2ang}

By the definition of the internal angle we have for $\phiint(0,C^{(p)}_w)$
\begin{eqnarray}
\phiint(0,C^{(p)}_w) &  = & \frac{1}{(2\pi)^{\frac{n}{2}}}\int_{-\1_{1\times (k-\keta)}\w_{\keta+1:k}-\1_{1\times n-k}|\w_{k+1:n}|\geq 0}e^{-\frac{\w^T\w}{2}}d\w\nonumber \\
&  = & \frac{1}{(2\pi)^{\frac{n-\keta}{2}}}\int_{-\1_{1\times (k-\keta)}\w_{\keta+1:k}-\1_{1\times n-k}|\w_{k+1:n}|\geq 0}e^{-\frac{\w_{\keta+1:n}^T\w_{\keta+1:n}}{2}}d\w.\nonumber \\
\label{eq:int2anal8}
\end{eqnarray}
Similarly to what we did above, one can again observe that the change $n\leftarrow n-\keta$ and $k\leftarrow k-\keta$ transforms the above angle to the corresponding one computed in \cite{Stojnicl1RegPosfinn}. Using the corresponding result from \cite{Stojnicl1RegPosfinn} and adjusting the dimensions we finally have
\begin{eqnarray}
\phiint(0,C^{(p)}_w)
& = & \frac{1}{2\pi}\lim_{x\rightarrow \infty}\lim_{\epsilon\rightarrow 0_+}(\int_{-\infty}^{-\epsilon} \lp1-i\erfi\lp\frac{t}{\sqrt{2}}\rp\rp^{n-k} e^{-\frac{lt^2}{2}}\frac{(1-e^{-itx})}{it}dt\nonumber \\
& & +\int_{\epsilon}^{\infty} \lp1-i\erfi\lp\frac{t}{\sqrt{2}}\rp\rp^{n-k} e^{-\frac{lt^2}{2}}\frac{(1-e^{-itx})}{it}dt).\label{eq:int2anal8a}
\end{eqnarray}

\subsection{External angles}
\label{sec:extang}

From (\ref{eq:intanal6}), we observe that there are two types of the external angles that we need to compute, $\phiext(F^{(l,p)}_1,C^{(p)}_w)$ and $\phiext(C^{(p)}_w,C^{(p)}_w)$. Clearly, $\phiext(C^{(p)}_w,C^{(p)}_w)=1$. Computing $\phiext(F^{(l,p)}_1,C^{(p)}_w)$ is a bit more involved and below we show how it can be done.

\subsubsection{Determining $\phiext(F^{(l,p)}_1,C^{(p)}_w)$}
\label{sec:ext1ang}

The computation that we show below can be sped up substantially. However, to ensure a completeness we show pretty much step by step how one can adapt the corresponding results from \cite{Stojnicl1RegPosfinn}. By the definition of the external angles, we have that at any given face they represent the content/fraction of the subspace taken by the positive hull of the outward normals to the hyperplanes that meet at the face. For face $F^{(l,p)}_1$ we then have that the corresponding positive hull is given as
\begin{equation}\label{eq:ext1anal1}
phull^{(l)}_{ext,1}\triangleq  -pos\lp\begin{bmatrix}
                                       \0_{\keta\times 1} \\
                                      -\1_{(l+1-\keta)\times 1} \\
                                      -\1^{(s,1)}_{(n-l-1)\times 1}
                                    \end{bmatrix},\begin{bmatrix}
                                       \0_{\keta\times 1} \\
                                      -\1_{(l+1-\keta)\times 1} \\
                                      -\1^{(s,2)}_{(n-l-1)\times 1}
                                    \end{bmatrix},\dots,\begin{bmatrix}
                                      \0_{\keta\times 1} \\
                                      -\1_{(l+1-\keta)\times 1} \\
                                      -\1^{(s,2^{n-l-1})}_{(n-l-1)\times 1}
                                    \end{bmatrix}\rp,\1^{(s,i)}_{(n-l-1)\times 1}\in \{-1,1\}^{n-l-1},
\end{equation}
and $\1^{(s,i)}_{(n-l-1)\times 1}\neq \1^{(s,j)}_{(n-l-1)\times 1}$ or any $i\neq j$ and $1\leq i,j\leq 2^{n-l+1}$. Due to symmetry one can work with $-phull^{(l)}_{ext,1}$ instead of $phull^{(l)}_{ext,1}$. Closely following \cite{Stojnicl1RegPosfinn} and relying on the ``Gaussian coordinates in an orthonormal basis" strategy we can choose the columns of the following matrix as a convenient orthonormal basis
\begin{equation}
B^{(l)}_{ext,1}=\begin{bmatrix}
e_{l+2} &  e_{l+3} & \dots & e_{n} & \begin{bmatrix}
                                      \0_{\keta\times 1} \\
                                       -\1_{(l+1-\keta)\times 1} \\
                                       \0_{(n-l-1)\times 1}
                                     \end{bmatrix}\frac{1}{\sqrt{l+1-\keta}}
          \end{bmatrix}.\label{eq:ext1anal2}
\end{equation}
It easily follows that $-phull^{(l)}_{ext,1}$ is indeed located in the subspace spanned by the columns of $B^{(l)}_{ext,1}$. Let $\g$ ($\g\in \mR^{n-l}$) be the coordinates in the basis determined by the columns of $B^{(l)}_{ext,1}$ and let $G^{(l)}_{ext,1}$ be such that
\begin{equation}
G^{(l)}_{ext,1} =\{\g\in \mR^{n-l}| B^{(l)}_{ext,1}\g\in -phull^{(l)}_{ext,1} \}.\label{eq:ext1anal3}
\end{equation}
To facilitate determining $G^{(l)}_{ext,1}$ we will introduce
\begin{equation}
D^{(l)}_{ext,1}=\begin{bmatrix}
\begin{bmatrix}
 \0_{\keta\times 1} \\
                                      -\1_{(l+1-\keta)\times 1} \\
                                      -\1^{(s,1)}_{(n-l-1)\times 1}
                                    \end{bmatrix} & \begin{bmatrix}
                                     \0_{\keta\times 1} \\
                                      -\1_{(l+1-\keta)\times 1} \\
                                      -\1^{(s,2)}_{(n-l-1)\times 1}
                                    \end{bmatrix} & \dots & \begin{bmatrix}
                                     \0_{\keta\times 1} \\
                                      -\1_{(l+1-\keta)\times 1} \\
                                      -\1^{(s,2^{n-l-1})}_{(n-l-1)\times 1}
                                    \end{bmatrix}\end{bmatrix},\label{eq:ext1anal4}
\end{equation}
and recognize that
\begin{eqnarray}
-phull^{(l)}_{ext,1} & = & \{D^{(l)}_{ext,1}\g^{(D)}|\g^{(D)}\geq 0,\g^{(D)}\in \mR^{2^{n-l-1}}\}\nonumber \\
& = & \{\begin{bmatrix}
\begin{bmatrix}
                                      -\1_{(l+1)\times 1} \\
                                      -\1^{(s,1)}_{(n-l-1)\times 1}
                                    \end{bmatrix} & \begin{bmatrix}
                                      -\1_{(l+1)\times 1} \\
                                      -\1^{(s,2)}_{(n-l-1)\times 1}
                                    \end{bmatrix} & \dots & \begin{bmatrix}
                                      -\1_{(l+1)\times 1} \\
                                      -\1^{(s,2^{n-l-1})}_{(n-l-1)\times 1}
                                    \end{bmatrix}\end{bmatrix}\g^{(D)}|\g^{(D)}\geq 0,\g^{(D)}\in \mR^{2^{n-l-1}}\}.\nonumber \\\label{eq:ext1anal5}
\end{eqnarray}
Combining (\ref{eq:ext1anal3}) and (\ref{eq:ext1anal5}) we obtain
\begin{equation}
G^{(l)}_{ext,1} =\{\g\in \mR^{n-l}| \exists \g^{(D)}\in\mR^{2^{n-l-1}}, \g^{(D)}\geq 0, \quad \mbox{and} \quad B^{(l)}_{ext,1}\g= D^{(l)}_{ext,1}\g^{(D)} \}.\label{eq:ext1anal6}
\end{equation}
Now, looking at the first $(l+1)$ equations in $B^{(l)}_{ext,1}\g= D^{(l)}_{ext,1}\g^{(D)}$ we also find that they imply
\begin{equation}
\g_{n-l}\frac{1}{\sqrt{l+1-\keta}}=\sum_{i=1}^{2^{n-l-1}}\g^{(D)}_{i}\geq 0.\label{eq:ext1anal7}
\end{equation}
Also, for $j\in\{2,3,\dots,n-l\}$, one can proceed line by line as in \cite{Stojnicl1RegPosfinn}, and, through analyzing $(l+j)$-th equations in $B^{(l)}_{ext,1}\g= D^{(l)}_{ext,1}\g^{(D)}$, obtain
\begin{equation}
G^{(l)}_{ext,1} =\{\g\in \mR^{n-l}| \g_{n-l}\geq 0, -\g_{n-l}\frac{1}{\sqrt{l+1-\keta}}\leq \g_{j-1}\leq \g_{n-l}\frac{1}{\sqrt{l+1-\keta}}, j\in\{2,3,\dots,n-l\} \}.\label{eq:ext1anal9}
\end{equation}
Moreover, following \cite{Stojnicl1RegPosfinn}, we have that (\ref{eq:ext1anal9}) is actually a complete characterization of $G^{(l)}_{ext,1}$, i.e. there are no other restrictions on $\g$. After switching to the Gaussian coordinates we finally have for $\phiext(F^{(l,p)}_1,C^{(p)}_w)$
\begin{multline}
\phiext(F^{(l,p)}_1,C^{(p)}_w)   = \frac{1}{(2\pi)^{\frac{n-l}{2}}}\int_{\g\in G^{(l)}_{ext,1}}e^{-\frac{\g^T\g}{2}}d\g  \\
 =  \frac{1}{(2\pi)^{\frac{l}{2}}}\int_{\g_{n-l}\geq 0}e^{-\frac{\g_{n-l}^2}{2}}\lp \prod_{j=2}^{n-l}\frac{1}{(2\pi)^{\frac{l}{2}}}\int_{ -\g_{n-l}\frac{1}{\sqrt{l+1-\keta}}}^{\g_{n-l}\frac{1}{\sqrt{l+1-\keta}}}e^{-\frac{\g_{j-1}^2}{2}}d\g_{j-1}\rp d\g_{n-l} \\
 =  \frac{1}{(2\pi)^{\frac{l}{2}}}\int_{\g_{n-l}\geq 0}e^{-\frac{\g_{n-l}^2}{2}}\lp \frac{1}{2}\erfc\lp\frac{-\g_{n-l}}{\sqrt{2}\sqrt{l+1-\keta}}\rp
-\frac{1}{2}\erfc\lp\frac{\g_{n-l}}{\sqrt{2}\sqrt{l+1-\keta}}\rp\rp^{n-l-1} d\g_{n-l}.
\label{eq:ext1anal10}
\end{multline}
The above calculations provide all the ingredients needed to compute $p^{(p)}_{err}(k,m,n,k_{\eta})$. The following theorem connects all the pieces together.
\begin{theorem}(Exact partial $\ell_1$'s performance characterization -- finite dimensions)
Let $A$ be an $m\times n$ matrix in (\ref{eq:system})
with i.i.d. standard normal components (or, alternatively, with the null-space uniformly distributed in the Grassmanian). Let
the unknown $\x$ in (\ref{eq:system}) be $k$-sparse. Further, let $supp(\x)$ and the signs of the nonzero components of $\x$ be arbitrarily chosen but fixed. Also, let $\Pi\subset supp(\x),|\Pi|=\keta$, and let $p^{(p)}_{err}(k,m,n,k_{\eta})$ be the probability that the solutions of (\ref{eq:l0}) and (\ref{eq:l1imp}) do \emph{not} coincide. Then
\begin{eqnarray}
p^{(p)}_{err}(k,m,n,k_{\eta}) & = & 2 ( \sum_{l=m+2j+1,j\in \mN_0}^{n-1} c^{(l)}_1\phiint(0,F^{(l,p)}_1)\phiext(F^{(l,p)}_1,C^{(p)}_w)+\phiint(0,C^{(p)}_w)\phiext(C^{(p)}_w,C^{(p)}_w)),\nonumber \\
\label{eq:finalthm1}
\end{eqnarray}
where $c^{(l)}_1$, $\phiint(0,F^{(l,p)}_1)$, $\phiint(0,C^{(p)}_w)$, and $\phiext(F^{(l,p)}_1,C^{(p)}_w)$ are as given in (\ref{eq:intanal4}), (\ref{eq:int1anal11}), (\ref{eq:int2anal8a}), and (\ref{eq:ext1anal10}),  respectively.
\label{thm:finalperr}
\end{theorem}
\begin{proof}
  Follows from the above discussion.
\end{proof}

\subsection{Simulations and theoretical results -- partial $\ell_1$}
\label{sec:posthnumresults}

\begin{figure}[htb]
\begin{minipage}[b]{.5\linewidth}
\centering
\centerline{\epsfig{figure=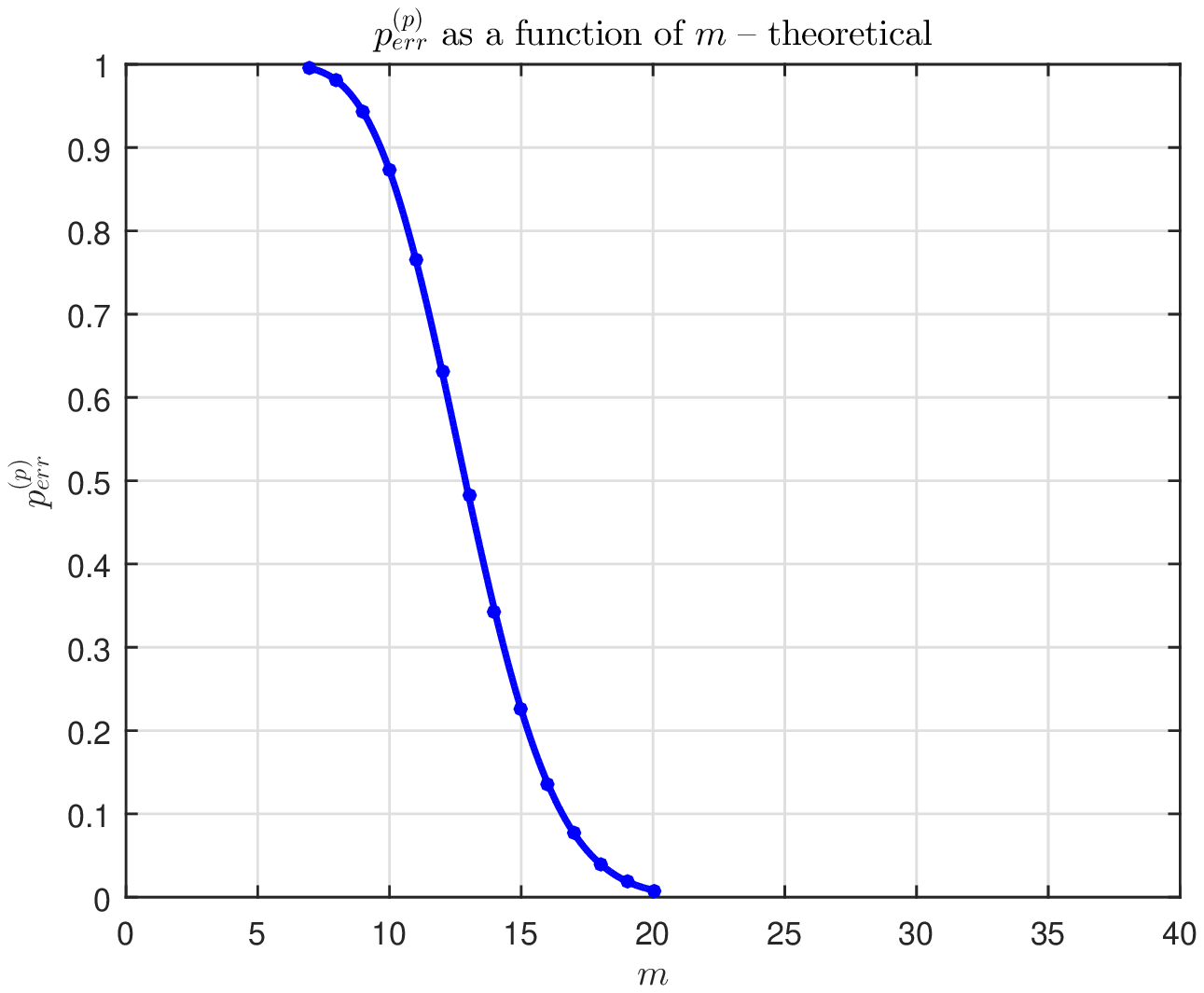,width=9cm,height=7cm}}
\end{minipage}
\begin{minipage}[b]{.5\linewidth}
\centering
\centerline{\epsfig{figure=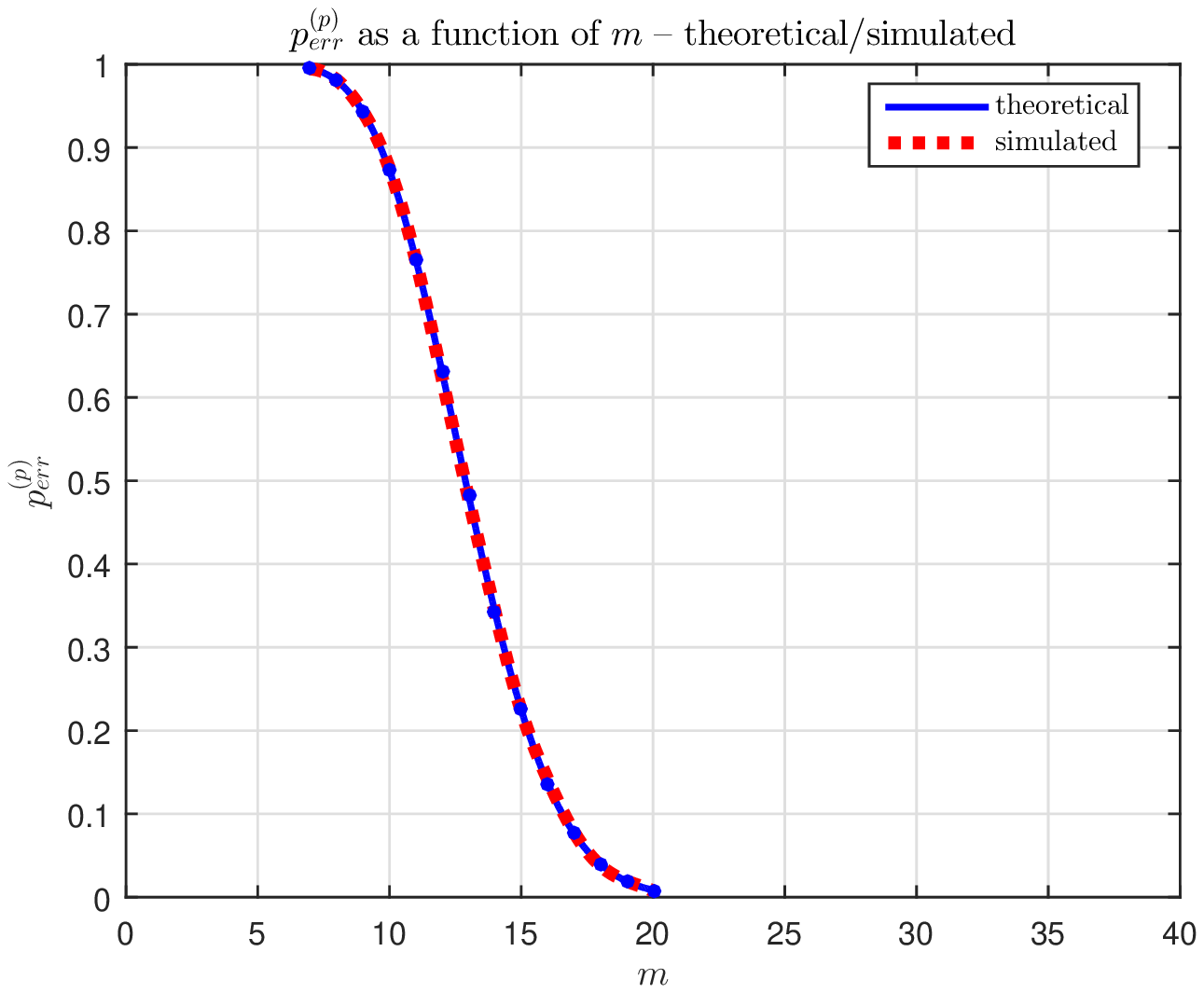,width=9cm,height=7cm}}
\end{minipage}
\caption{$p^{(p)}_{err}(k,m,n,k_{\eta})$ as a function of $m$ ($k=6$, $n=40$, $k_\eta=3$); left -- theory; right -- simulations}
\label{fig:l1regnonperr}
\end{figure}

In this section we will determine the concrete values for $p^{(p)}_{err}(k,m,n,k_{\eta})$ based on what is proven in Theorem \ref{thm:finalperr}. In Figure \ref{fig:l1regnonperr} we show both, the simulated and the theoretical values for $p^{(p)}_{err}(k,m,n,k_{\eta})$ (the theoretical values are, of course, obtained based on Theorem \ref{thm:finalperr}). We fixed $k=6$, $n=40$, and $\keta=3$ and varied/increased $m$ so that $p^{(p)}_{err}(k,m,n,k_{\eta})$ changes from one to zero. Figure \ref{fig:l1regnonperr} is complemented by Table \ref{tab:l1regnonperrtab1} where we show the numerical values for $p^{(p)}_{err}(k,m,n,k_{\eta})$ (again, both, simulated and theoretical) obtained for several concrete values of quadruplets $(k,m,n,k_{\eta})$ (we also show the number of numerical experiments that were run as well as the number of them that did not result in having the solution of (\ref{eq:l1imp}) match the a priori known to be nonnegative solution of (\ref{eq:l0})). We observe an excellent agreement between the simulated and theoretical results.

\begin{table}[h]
\caption{Simulated and theoretical results for $p^{(p)}_{err}(k,m,n,k_{\eta})$; $k=6$, $n=40$, $k_\eta=3$}\vspace{.1in}
\hspace{-0in}\centering
\begin{tabular}{||c||c|c|c|c|c|c||}\hline\hline
$m$                    & $ 10 $ & $ 11 $ & $ 12 $ & $ 13 $ & $ 14 $ & $ 15 $ \\ \hline \hline
$\#$ of failures       & $ 9929 $ & $ 4227 $ & $ 6235 $ & $ 6332 $ & $ 4263 $ & $ 2671 $ \\ \hline
$\#$ of repetitions    & $ 11354 $ & $ 5528 $ & $ 10010 $ & $ 13217 $ & $ 12427 $ & $ 11934 $ \\ \hline \hline
$p^{(p)}_{err}(k,m,n,k_{\eta})$ -- simulation& $ \bl{\mathbf{0.8745}} $ & $ \bl{\mathbf{0.7647}} $ & $ \bl{\mathbf{0.6229}} $ & $ \bl{\mathbf{0.4791}} $ & $ \bl{\mathbf{0.3430}} $ & $ \bl{\mathbf{0.2238}} $ \\ \hline \hline
$p^{(p)}_{err}(k,m,n,k_{\eta})$ -- theory    & $ \mathbf{0.8722} $ & $ \mathbf{0.7652} $ & $ \mathbf{0.6300} $ & $ \mathbf{0.4826} $ & $ \mathbf{0.3429} $ & $ \mathbf{0.2251} $ \\ \hline \hline
\end{tabular}
\label{tab:l1regnonperrtab1}
\end{table}

\subsection{Asymptotics}
\label{sec:asym}

In our companion paper \cite{Stojnicl1HidParasymldp} we look at the partial $\ell_1$ in an asymptotic scenario. In such a scenario, the most typical is the so-called linear regime where the systems dimensions grow larger in a linearly proportional fashion. As mentioned earlier, in an asymptotic linear regime the standard and positive $\ell_1$ exhibit the so-called phase-transition phenomenon that was fully settled through \cite{DonohoPol,DonohoUnsigned,StojnicCSetam09,StojnicUpper10}.  Similar behavior turns out to be present in the case of the partial $\ell_1$ as well. In precisely this same linear asymptotic regime \cite{StojnicICASSP10knownsupp,StojnicTowBettCompSens13} settled the phase-transition behavior of the partial $\ell_1$. The companion paper, \cite{Stojnicl1HidParasymldp}, then went much further and determined the partial $\ell_1$'s LDP behavior. It did so by developing two fundamentally different mathematical approaches, one that is purely probabilistic in nature and another one that relies on some of the considerations presented here. The analysis presented in \cite{Stojnicl1HidParasymldp} is somewhat technical and a detailed discussion about it is beyond our interests here. However, below we sketch how the finite dimensional characterizations transform into the asymptotic ones eventually utilized in \cite{Stojnicl1HidParasymldp}.

Before proceeding with the main ideas we quickly recall what the key difference between the PT and the LDP is. Namely, settling the PT phenomenon usually means determining the PT curve in $(\alpha,\beta)$ plane. There are several equivalent ways how one can look at that problem. For example, one of them would be finding for any $\beta\in(0,1)$ a critical (the smallest) $\alpha$ such that for any $\alpha$ above it the algorithm (in our case here, the partial $\ell_1$) succeeds with overwhelming probability. Settling the LDP phenomenon, on the other hand, assumes not only determining the critical $\alpha$ but also precisely characterizing the rate at which probability that the partial $\ell_1$ fails or succeeds goes to zero as one moves around the critical $\alpha$ in the so-called transition zone. In other words, it assumes determination of $I^{(p)}_{ldp}(\alpha,\beta;\eta)=\lim_{n\rightarrow\infty}\frac{\log(p^{(p)}_{err}(k,m,n,k_{\eta}))}{n}$, where say $\keta=\eta k$. In \cite{Stojnicl1HidParasymldp} we settle that problem. The analysis that we presented in earlier sections can be used as a starting for what is presented in \cite{Stojnicl1HidParasymldp} (below we provide a brief sketch as to how to bridge between the two, the above analysis and \cite{Stojnicl1HidParasymldp}; it essentially boils down to
transforming the above $p^{(p)}_{err}(k,m,n,k_{\eta})$ to a more convenient characterization in an infinite dimensional setting).

In what follows we will be dealing with the linear regime. To ensure that everything is properly scaled we will assume that $k=\beta n, m=\alpha n, \keta=\eta\beta n$, and $l=\rho n$ (where $\beta$, $\alpha$, $\eta$,and $\rho$ are fixed constants independent of $n$). Following the strategy presented in \cite{Stojnicl1RegPosfinn}, we, as $n\rightarrow\infty$, from (\ref{eq:finalthm1}) have
\begin{equation}
\lim_{n\rightarrow \infty}\frac{\log(p^{(p)}_{err}(k,m,n,k_{\eta}))}{n}  =  \max\{\max_{\rho\geq \alpha}\lim_{n\rightarrow \infty} \frac{\log(\zeta^{(\infty)}_1)}{n},\lim_{n\rightarrow \infty}\frac{\log(\zeta^{(\infty)}_2)}{n}\},
\label{eq:asym1}
\end{equation}
where
\begin{eqnarray} \label{eq:asym2}
  \lim_{n\rightarrow \infty} \frac{\log(\zeta^{(\infty)}_1)}{n} &=&
  \lim_{n\rightarrow \infty} \lp\frac{\log(c^{(l)}_{1})}{n}+
  \frac{\log(\phiint(0,F^{(l,p)}_1))}{n}+\frac{\log(\phiext(F^{(l,p)}_1,C^{(p)}_w))}{n}\rp \nonumber \\
  \lim_{n\rightarrow \infty} \frac{\log(\zeta^{(\infty)}_2)}{n} &=&
  \lim_{n\rightarrow \infty} \lp
  \frac{\log(\phiint(0,C^{(p)}_w))}{n}+\frac{\log(\phiext(C^{(p)}_w,C^{(p)}_w))}{n}\rp.
\end{eqnarray}
From (\ref{eq:intanal4}) we have
\begin{eqnarray}
\lim_{n\rightarrow \infty} \frac{\log(c^{(l)}_{1})}{n}=\lim_{n\rightarrow \infty}\frac{\log\lp 2^{l-k+1}\binom{n-k}{n-l-1}\rp}{n}=-(1-\beta)H\lp\frac{1-\rho}{1-\beta}\rp +(\rho-\beta)\log(2).\label{eq:asym3}
\end{eqnarray}
We recall on $F^{(l,p)}_1$ and introduce $F^{(l,p)}_2$
\begin{eqnarray}\label{eq:asym5a}
F^{(l,p)}_1 & \triangleq & F^{(l,p)}_1(k,\keta)= \{\w\in \mR^n|-\sum_{i=\keta+1}^k \w_i = \sum_{i=k+1}^{l+1}\w_{i},\w_{k+1:l+1}\geq 0,\w_{l+2:n}= 0\} \nonumber \\
F^{(l,p)}_2 & \triangleq & F^{(l,p)}_2(k,\keta)= \{\w\in \mR^n|-\sum_{i=\keta+1}^k \w_i\geq \sum_{i=k+1}^{l}\w_{i},\w_{k+1:l}\geq 0,\w_{l+1:n}= 0\}.
\end{eqnarray}
From \cite{Stojnicl1RegPosfinn} we also recall
\begin{eqnarray}\label{eq:asym5b}
F^{(l,+)}_1 & \triangleq & F^{(l,+)}_1(k)= \{\w\in \mR^n|-\sum_{i=1}^k \w_i = \sum_{i=k+1}^{l+1}\w_{i},\w_{k+1:l+1}\geq 0,\w_{l+2:n}= 0\} \nonumber \\
F^{(l,+)}_2 & \triangleq & F^{(l,+)}_2(k)= \{\w\in \mR^n|-\sum_{i=1}^k \w_i\geq \sum_{i=k+1}^{l}\w_{i},\w_{k+1:l}\geq 0,\w_{l+1:n}= 0\}.
\end{eqnarray}
Moreover in \cite{Stojnicl1RegPosfinn} it was also shown that
\begin{eqnarray}
\lim_{n\rightarrow \infty} \frac{\log(\phiint(0,F^{(l,+)}_1))}{n}  =  \lim_{n\rightarrow \infty} \frac{\log(\phiint(0,F^{(l,+)}_2))}{n},
\label{eq:asym5c}
\end{eqnarray}
which by (\ref{eq:asym5b}) is equivalent to
\begin{eqnarray}
\lim_{n\rightarrow \infty} \frac{\log(\phiint(0,F^{(l,+)}_1(k)))}{n}  =  \lim_{n\rightarrow \infty} \frac{\log(\phiint(0,F^{(l,+)}_2(k)))}{n}.
\label{eq:asym5d}
\end{eqnarray}
Combining (\ref{eq:asym5a}), (\ref{eq:asym5b}), (\ref{eq:asym5c}), (\ref{eq:asym5d}), (\ref{eq:posint1anal8}), and considerations from \cite{Stojnicl1RegPosfinn} one obtains
\begin{eqnarray}
\lim_{n\rightarrow \infty} \frac{\log(\phiint(0,F^{(l,p)}_1))}{n}  & = &  \lim_{n\rightarrow \infty} \frac{\log(\phiint(0,F^{(l,p)}_1(k,\keta))}{n} \nonumber \\
& = &
\lim_{n\rightarrow \infty} \frac{\log(\phiint(0,F^{(l-\keta,+)}_1(k-\keta)))}{n}\nonumber \\
& = &  \lim_{n\rightarrow \infty} \frac{\log(\phiint(0,F^{(l-\keta,+)}_2(k-\keta)))}{n}\nonumber \\.
& = & \lim_{n\rightarrow \infty} \frac{\log(\phiint(0,F^{(l,p)}_2(k,\keta)))}{n}\nonumber \\
& = & \lim_{n\rightarrow \infty} \frac{\log(\phiint(0,F^{(l,p)}_2))}{n}.
\label{eq:asym5e}
\end{eqnarray}
We now also introduce, as in \cite{Stojnicl1RegPosfinn},  a mathematical object, $\phiint(0,F^{(n)}_1)$ ($F^{(n)}_1$ is clearly not a face of $C^{(p)}_w$), in the following way
\begin{equation}
\phiint(0,F^{(n)}_1)   \triangleq
\frac{\sqrt{n+1-\keta}}{(2\pi)^{\frac{n-\keta}{2}}}\int_{-\1_{1\times (n-\keta)}\w_{\keta+1:n}\geq 0,\w_{k+1:n}\geq 0}e^{-\frac{\w_{\keta+1:n}^T\w_{\keta+1:n}+\w_{\keta+1:n}^T( -\1_{1\times n})^T
             (-\1_{1\times n})\w_{\keta+1:n}}{2}}d\w_{\keta+1:n}.\\
\label{eq:asym4b}
\end{equation}
Using the same line of reasoning as in \cite{Stojnicl1RegPosfinn} we then have
\begin{eqnarray}
\lim_{n\rightarrow \infty} \frac{\log(\phiint(0,F^{(n-1,p)}_1))}{n} &  = & \lim_{n\rightarrow \infty} \frac{\log(\phiint(0,F^{(n,p)}_1))}{n}\nonumber \\
&  = & \lim_{n\rightarrow \infty} \frac{\log(\phiint(0,F^{(n,p)}_2))}{n}\nonumber \\
&  =  & \lim_{n\rightarrow \infty} \frac{\log\lp\frac{1}{(2\pi)^{\frac{n-\keta}{2}}}\int_{-\1_{1\times (n-\keta)}\w_{\keta+1:n}\geq 0,\w_{k+1:n}\geq 0}e^{-\frac{\w_{\keta+1:n}^T\w_{\keta+1:n}}{2}}d\w_{\keta+1:n}\rp}{n}\nonumber \\
&  =  & \lim_{n\rightarrow \infty} \frac{\log\lp\frac{1}{2^{n-k}}\phiint(0,C^{(p)}_w)\rp}{n},\nonumber \\\label{eq:asym4c}
\end{eqnarray}
From (\ref{eq:asym1}), (\ref{eq:asym2}), (\ref{eq:asym3}), and (\ref{eq:asym4c}) we obtain
\begin{equation}
\lim_{n\rightarrow \infty}\frac{\log(p^{(hp)}_{err}(k,m,n,k_{\eta}))}{n}  =  \max_{\rho\geq \alpha}\lim_{n\rightarrow \infty} \frac{\log(\zeta^{(\infty)}_1)}{n}.
\label{eq:asym5}
\end{equation}
From (\ref{eq:asym5e}) one also has
\begin{eqnarray}
\lim_{n\rightarrow \infty} \frac{\log(\phiint(0,F^{(l,p)}_1))}{n} &  = & \lim_{n\rightarrow \infty} \frac{\log(\phiint(0,F^{(l,p)}_2))}{n} \nonumber \\
& = & \lim_{n\rightarrow \infty} \frac{\log\lp\frac{1}{(2\pi)^{\frac{l-\keta}{2}}}\int_{-\1_{1\times (l-\keta)}\w_{\keta+1:l}\geq 0,\w_{k+1:l}\geq 0}e^{-\frac{\w_{\keta+1:l}^T\w_{\keta+1:l}}{2}}d\w_{\keta+1:l}\rp}{n}\nonumber \\
& = & \lim_{n\rightarrow \infty} \frac{\log\lp\frac{2^{l-k}}{2^{l-k}(2\pi)^{\frac{l-\keta}{2}}}\int_{-\1_{1\times (l-\keta)}\w_{\keta+1:l}\geq 0,\w_{k+1:l}\geq 0}e^{-\frac{\w_{\keta+1:l}^T\w_{\keta+1:l}}{2}}d\w_{\keta+1:l}\rp}{n}\nonumber \\
&  = & \lim_{n\rightarrow \infty} \frac{\log\lp\frac{1}{2^{l-k}}P(-\1_{1\times (l-\keta)}\w_{\keta+1:l}\geq 0)\rp}{n},\label{eq:posasym6}
\end{eqnarray}
where on the right side of the last equality one can think of the elements of $\w_{\keta+1:k}$ as being the i.i.d. standard normals and the elements of $\w_{k+1:l}$ as being the i.i.d. standard half normals. One can now continue following what was presented in \cite{Stojnicl1RegPosfinn} to obtain
\begin{eqnarray}
\lim_{n\rightarrow \infty} \frac{\log(\phiint(0,F^{(l,p)}_1))}{n}
&  = & \lim_{n\rightarrow \infty} \frac{\log\lp\frac{1}{2^{l-k}}P(-\1_{1\times (l-\keta)}\w_{\keta+1:l}\geq 0)\rp}{n}\nonumber \\
& = &  \min_{\mu_y\geq 0} \lim_{n\rightarrow \infty} \frac{\log(\mE e^{-\mu_y\1_{1\times (l-\keta)}\w_{\keta+1:l}})}{n}-(\rho-\beta)\log(2)\nonumber \\
& = &  \min_{\mu_y\geq 0} \lp(\rho-\beta) \log\lp\mE e^{-\mu_y\w_{k+1}}\rp+(1-\eta)\beta\frac{\mu_y^2}{2}\rp-(\rho-\beta)\log(2) \nonumber \\
& = &  \min_{\mu_y\geq 0} \lp(\rho-\beta) \log\lp\frac{2}{\sqrt{2\pi}}\int_{0}^{\infty} e^{-\frac{\w_{k+1}^2}{2}-\mu_y\w_{k+1}}d\w_{k+1}\rp+(1-\eta)\beta\frac{\mu_y^2}{2}\rp\nonumber \\
& & -(\rho-\beta)\log(2) \nonumber \\
& = &  \min_{\mu_y\geq 0} \lp(\rho-\beta) \log\lp\erfc\lp\frac{\mu_y}{\sqrt{2}}\rp\rp+(\rho-\eta\beta)\frac{\mu_y^2}{2}\rp-(\rho-\beta)\log(2) \nonumber \\
& = &  \min_{\mu_y\geq 0} \lp(\rho-\beta) \log(\erfc(\mu_y))+(\rho-\eta\beta)\mu_y^2\rp-(\rho-\beta)\log(2).\label{eq:asym7}
\end{eqnarray}
From (\ref{eq:ext1anal10}) we obtain
\begin{eqnarray}
\lim_{n\rightarrow \infty} \frac{\log(\phiext(F^{(l,p)}_1,C^{(p)}_w))}{n}
& = & \max_{\g_{n-l}\geq 0} \lp -\frac{\g_{n-l}^2}{2}+(1-\rho)\log\lp \frac{1}{2}\erfc\lp\frac{-\g_{n-l}}{\sqrt{2}\sqrt{\rho-\eta\beta}}\rp
-\frac{1}{2}\erfc\lp\frac{\g_{n-l}}{\sqrt{2}\sqrt{\rho-\eta\beta}}\rp\rp\rp \nonumber \\
& = & \max_{\g_{n-l}\geq 0} \lp -(\rho-\eta\beta) \g_{n-l}^2+(1-\rho)\log\lp \frac{1}{2}\erfc(-\g_{n-l})-\frac{1}{2}\erfc(\g_{n-l})\rp\rp,
\label{eq:asym8}
\end{eqnarray}
Finally, a combination of (\ref{eq:int1anal11}), (\ref{eq:ext1anal10}), (\ref{eq:asym2}), (\ref{eq:asym3}), (\ref{eq:asym7}), and (\ref{eq:asym8}) transforms (\ref{eq:asym5}) into the following
\begin{eqnarray}
\lim_{n\rightarrow \infty}\frac{\log(p^{(p)}_{err}(k,m,n,k_{\eta}))}{n}
 & = & \max_{\rho\geq \alpha}\lim_{n\rightarrow \infty} \frac{\log(\zeta^{(\infty)}_1)}{n}\nonumber \\
 & = & \max_{\rho\geq \alpha} (-(1-\beta)H\lp\frac{1-\rho}{1-\beta}\rp+(\rho-\beta)\log(2) \nonumber \\
 & & + \min_{\mu_y\geq 0} \lp(\rho-\beta) \log(\erfc(\mu_y))+(\rho-\eta\beta)\mu_y^2\rp-(\rho-\beta)\log(2)\nonumber \\.
& & + \max_{\g_{n-l}\geq 0} \lp -(\rho-\eta\beta) \g_{n-l}^2+(1-\rho)\log\lp \frac{1}{2}\erfc(-\g_{n-l})-\frac{1}{2}\erfc(\g_{n-l})\rp\rp). \nonumber \\
\label{eq:asym12}
\end{eqnarray}
Now, consider a given $\beta$, and let, for such a $\beta$, $\alpha_w$ be the $\alpha$ that produces $\lim_{n\rightarrow \infty}\frac{\log(p^{(p)}_{err}(k,m,n,k_{\eta}))}{n} =0$ (such an $\alpha$ always exists; although this is rather obvious, it is proven rigorously in \cite{Stojnicl1HidParasymldp}). Following the reasoning from \cite{Stojnicl1RegPosfinn} one can conclude that if $\alpha\geq \alpha_w$ then $\rho=\alpha$ is optimal in (\ref{eq:asym12}) and  the optimization over $\rho$ in (\ref{eq:asym12}) can be removed. On can also apply the same reasoning even if $\alpha\leq \alpha_w$. The only difference is that in such a scenario one focuses on the complementary version of (\ref{eq:anal3}), i.e. one focuses on
\begin{equation}
P(G^{(sub)}\cap C^{(p)}_w\neq \emptyset)=1-2\sum_{l=m-2j-1,j\in \mN_0,l\geq k-1} \sum_{F^{(l)}\in \calF^{(l)}}\phiint(0,F^{(l)})\phiext(F^{(l)},C^{(p)}_w)=1-p^{(p)}_{cor},\label{eq:asym13}
\end{equation}
where $p^{(p)}_{cor}$ is the probability of being correct, i.e. the probability that the solution of (\ref{eq:l1imp}) is the $k$-sparse solution of (\ref{eq:l0}) and its decay rate is given by (\ref{eq:asym12}) with $\rho=\alpha$. (\ref{eq:asym12}) is then sufficient to fully determine numerically PT and LDP curves of the standard $\ell_1$. What \cite{Stojnicl1HidParasymldp} does is, however, way beyond that; namely, \cite{Stojnicl1HidParasymldp}, among other things, explicitly analytically solves (\ref{eq:asym12}).

\section{Hidden partial $\ell_1$ -- finite dimensional analysis}
\label{sec:hidl1}

In this section we analyze the hidden partial $\ell_1$ from (\ref{eq:l1imphidden}). As we will soon see, quite a lot of what was done in Section \ref{sec:posl1} can be reused here. Before doing that we need to establish a couple of facts analogous to those that Section \ref{sec:posl1} relies on. We will again, for the concreteness of the exposition and without loss of generality, assume that the elements $\x_{k+1},\x_{k+2},\dots,\x_{n}$ of $\x$ are equal to zero and that the elements $\x_{1},\x_{2},\dots,\x_k$ have fixed signs, say all positive. Also, for concreteness and without loss of generality, let $\kappa=\{1,2,\dots,\keta,n-(k-\keta)+1,,n-(k-\keta)+2,\dots,n\}$. The following is a hidden partial $\ell_1$ adaptation of the result proven for the general $\ell_1$ in \cite{StojnicCSetam09,StojnicUpper10,StojnicICASSP09} (essentially an analogue to Theorem \ref{thm:posthmregposcond}).
\begin{theorem}(\cite{StojnicTowBettCompSens13} Nonzero elements of $\x$ have fixed signs and location)
Assume that an $m\times n$ system matrix $A$ is given. Let $\x$
be a $k$ sparse vector. Also let $\x_{k+1}=\x_{k+2}=\dots=\x_{n}=0$. Let the signs of $\x_{1},\x_{2},\dots,\x_k$ be fixed, say all positive. Also, let $\kappa=\{1,2,\dots,\keta,n-(k-\keta)+1,,n-(k-\keta)+2,\dots,n\}$. Further, assume that $\y=A\x$ and that $\w$ is
a $n\times 1$ vector. If
\begin{equation}
(\forall \w\in \textbf{R}^n | A\w=0) \quad  -\sum_{i=\keta+1}^{k} \w_i<\sum_{i=k+1}^{n-(k-\keta)}|\w_{i}|,
\label{eq:hidthmcond1}
\end{equation}
then the solutions of (\ref{eq:l0}) and (\ref{eq:l1imphidden}) coincide. Moreover, if
\begin{equation}
(\exists \w\in \textbf{R}^n | A\w=0) \quad  -\sum_{i=\keta+1}^{k} \w_i\geq \sum_{i=k+1}^{n-(k-\keta)}|\w_{i}|,
\label{eq:hidthmcond2}
\end{equation}
then there is an $\x$ from the above set of $\x$'s with fixed location of nonzero elements such that the solution of (\ref{eq:l0}) is not the solution of (\ref{eq:l1imphidden}).
\label{thm:hidthmregposcond}
\end{theorem}
Analogously to defining $C^{(p)}_w$ in Section \ref{sec:posl1} we can now define
\begin{equation}
C^{(hp)}_w\triangleq C^{(hp)}_w(k,m,n,\keta) =\{\w\in \mR^n| \quad -\sum_{i=\keta+1}^k \w_i\geq \sum_{i=k+1}^{n-(k-\keta)}|\w_{i}|\}.\label{eq:hiddefSw}
\end{equation}
Now, comparing $C^{(hp)}$ to $C^{(p)}$, one can observe that
\begin{equation}
C^{(hp)}_w(k,m,n,\keta) = C^{(p)}_w(2k-\keta,m,n,k).\label{eq:hiddefSw}
\end{equation}
Following what was done in Section \ref{sec:posl1} and utilizing (\ref{eq:hiddefSw}) we have
\begin{eqnarray}
p^{(hp)}_{err}(k,m,n,k_{\eta}) & = & P(G^{(sub)}\cap C^{(hp)}_w(k,m,n,\keta)\neq \emptyset) \nonumber \\
& = & P(G^{(sub)}\cap C^{(p)}_w(2k-\keta,m,n,k)\neq \emptyset) \nonumber \\
& = & p^{(p)}_{err}(2k-\keta,m,n,k),\label{eq:hidanal3}
\end{eqnarray}
where $p^{(hp)}_{err}(k,m,n,k_{\eta})$ is the probability of error that the $k$-sparse solution of (\ref{eq:l0}) is not the solution of (\ref{eq:l1imphidden}). The following is then the hidden partial analogue to Theorem \ref{thm:finalperr}.
\begin{theorem}(Exact hidden partial $\ell_1$'s performance characterization -- finite dimensions)
Let $A$ be an $m\times n$ matrix in (\ref{eq:system})
with i.i.d. standard normal components (or, alternatively, with the null-space uniformly distributed in the Grassmanian). Let
the unknown $\x$ in (\ref{eq:system}) be $k$-sparse. Further, let $supp(\x)$ and the signs of the nonzero components of $\x$ be arbitrarily chosen but fixed. Also, let $\kappa$ be a set such that $|\kappa\cap supp(\x)|=\keta,|\kappa|=k$, and let $p^{(hp)}_{err}(k,m,n,k_{\eta})$ be the probability that the solutions of (\ref{eq:l0}) and (\ref{eq:l1imphidden}) do \emph{not} coincide. Then
\begin{equation}
p^{(hp)}_{err}(k,m,n,k_{\eta}) =  p^{(p)}_{err}(2k-\keta,m,n,k),
\label{eq:hidfinalthm1}
\end{equation}
where
$p^{(p)}_{err}(k,m,n,\keta)$ is as in Theorem \ref{thm:finalperr}.
\label{thm:hidfinalperr}
\end{theorem}
\begin{proof}
  Follows from the above discussion.
\end{proof}

\subsection{Simulations and theoretical results -- hidden partial $\ell_1$}
\label{sec:hidthnumresults}

\begin{figure}[htb]
\begin{minipage}[b]{.5\linewidth}
\centering
\centerline{\epsfig{figure=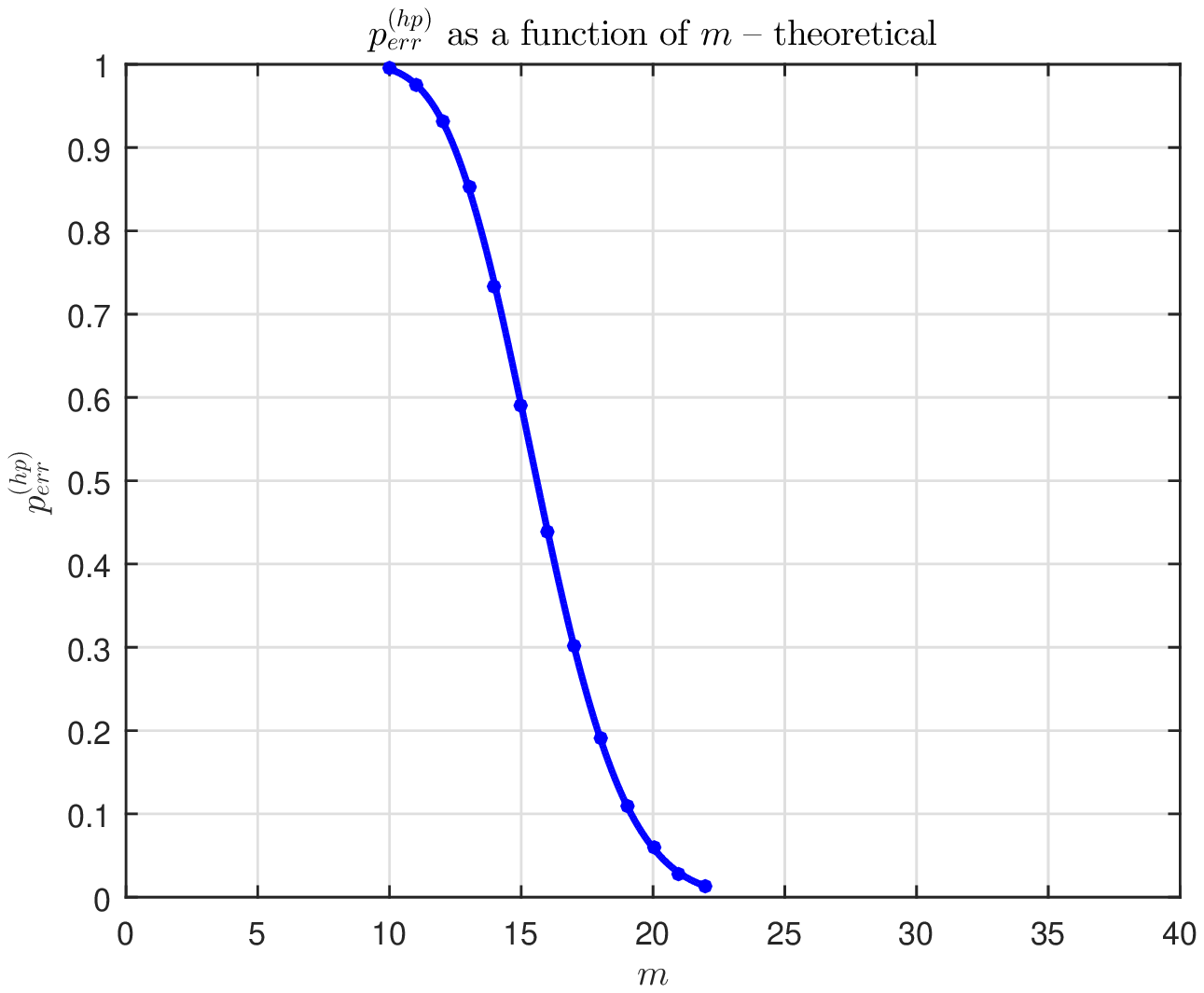,width=9cm,height=7cm}}
\end{minipage}
\begin{minipage}[b]{.5\linewidth}
\centering
\centerline{\epsfig{figure=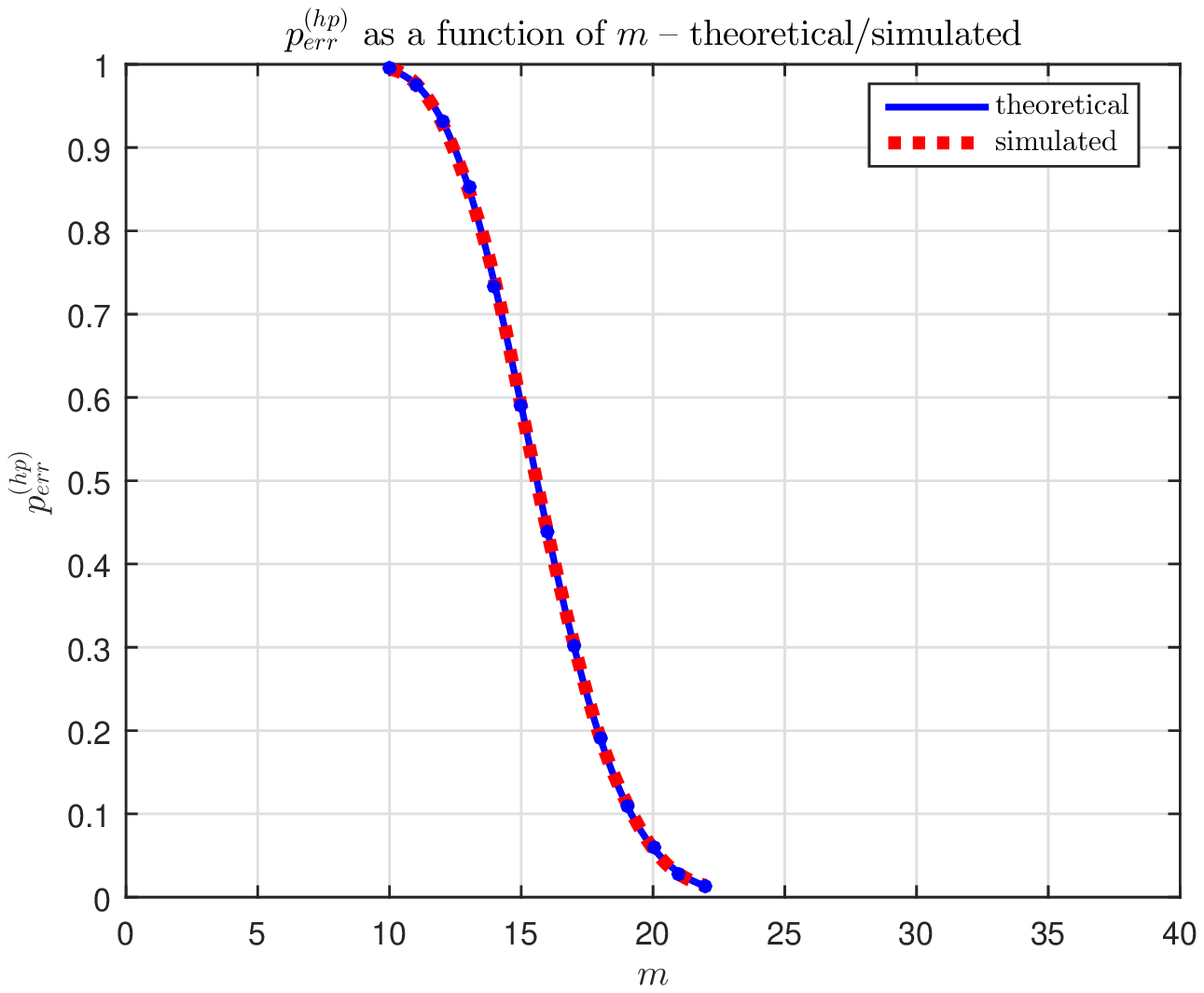,width=9cm,height=7cm}}
\end{minipage}
\caption{$p^{(hp)}_{err}(k,m,n,k_{\eta})$ as a function of $m$ ($k=6$, $n=40$, $k_\eta=3$); left -- theory; right -- simulations}
\label{fig:l1hidperr}
\end{figure}

Using Theorems \ref{thm:finalperr} and \ref{thm:hidfinalperr} one can obtain concrete theoretical values for $p^{(hp)}_{err}(k,m,n,k_{\eta})$. These values together with the corresponding ones that can be obtained through numerical simulations are presented in this section. First, in Figure \ref{fig:l1hidperr} we show the simulated and the theoretical plots of $p^{(hp)}_{err}(k,m,n,k_{\eta})$, and then in Table \ref{tab:l1hidperrtab1} we show the corresponding numerical values for several concrete quadruplets $(k,m,n,k_{\eta})$. To make presented results easier to analyze, in all experiments and theoretical calculations, we fixed $k=6$, $n=40$, and $\keta=3$ and varied/increased $m$ so that $p^{(hp)}_{err}(k,m,n,k_{\eta})$ changes from one to zero. Additionally, for each of the chosen quadruplets we also show the number of the numerical experiments that were run as well as the number of them that did not result in having the solution of (\ref{eq:l1imphidden}) match the $k$-sparse solution of (\ref{eq:l0}). As in Section \ref{sec:posl1}, we observe an excellent agreement between the simulated and the theoretical results.

\begin{table}[h]
\caption{Simulated and theoretical results for $p^{(hp)}_{err}(k,m,n,k_{\eta})$; $k=6$, $n=40$, $k_\eta=3$}\vspace{.1in}
\hspace{-0in}\centering
\begin{tabular}{||c||c|c|c|c|c|c||}\hline\hline
$m$                    & $ 13 $ & $ 14 $ & $ 15 $ & $ 16 $ & $ 17 $ & $ 18 $ \\ \hline \hline
$\#$ of failures       & $ 17154 $ & $ 11906 $ & $ 10621 $ & $ 7683 $ & $ 6036 $ & $ 3585 $ \\ \hline
$\#$ of repetitions    & $ 20203 $ & $ 16094 $ & $ 18036 $ & $ 17543 $ & $ 19857 $ & $ 18750 $ \\ \hline \hline
$p^{(hp)}_{err}(k,m,n,k_{\eta})$ -- simulation& $ \bl{\mathbf{0.8491}} $ & $ \bl{\mathbf{0.7398}} $ & $ \bl{\mathbf{0.5889}} $ & $ \bl{\mathbf{0.4380}} $ & $ \bl{\mathbf{0.3040}} $ & $ \bl{\mathbf{0.1912}} $ \\ \hline \hline
$p^{(hp)}_{err}(k,m,n,k_{\eta})$ -- theory    & $ \mathbf{0.8519} $ & $ \mathbf{0.7344} $ & $ \mathbf{0.5902} $ & $ \mathbf{0.4394} $ & $ \mathbf{0.3014} $ & $ \mathbf{0.1906} $ \\ \hline \hline
\end{tabular}
\label{tab:l1hidperrtab1}
\end{table}

\subsection{Asymptotics}
\label{sec:hidasym}

Following Section \ref{sec:posl1} and utilizing Theorem \ref{thm:hidfinalperr} we also easily find
\begin{equation}\label{eq:hidasym1}
I^{(hp)}_{ldp}(\alpha,\beta;\eta)=\lim_{n\rightarrow\infty}\frac{\log(p^{(hp)}_{err}(k,m,n,k_{\eta}))}{n}
=\lim_{n\rightarrow\infty}\frac{\log(p^{(p)}_{err}(2k-\keta,m,n,k))}{n}=I^{(p)}_{ldp}(\alpha,(2-\eta)\beta;(2-\eta)^{-1}),
\end{equation}
which is enough to bridge between the finite dimensional scenario considered here and the asymptotic one considered in \cite{Stojnicl1HidParasymldp}.

\section{Conclusion}
\label{sec:conc}

In this paper we studied finite dimensional random linear systems with sparse solutions. In particular, we focused on a couple of modifications of the standard $\ell_1$ heuristic and their performance analysis. The modifications that we considered are typical for scenarios where one has a bit of feedback as to what the unknown vectors are. As is well known, in linear systems with sparse solutions the key to solving the problem is determination of the unknown vector's support, i.e. the location of its nonzero components. In such problems a fairly useful feedback would be any knowledge about the true support. Here we first analyzed the scenario where a portion of the support is a priori known. To do that we utilized the well known partial $\ell_1$ modification of the standard $\ell_1$. We then considered a bit more practical scenario where a portion of the true support is known but it is so to say hidden in a larger set that is available. To handle that we introduced in \cite{StojnicTowBettCompSens13} the so-called hidden partial $\ell_1$ and here we provided its a finite dimensional performance analysis as well. The results obtained through the analyses show very precisely what kind of gains one can expect to see from the available feedback.

Both of these modifications (especially the hidden partial one) also seem as a promising tool in designing algorithms that could in certain scenarios outperform the standard $\ell_1$. When viewed this way, they in fact help transforming the original sparse solutions recovery problem into the partial support recovery. In that direction, the results provided here show precisely what kind of partial support recovery problem one should be able to solve and with what statistical guarantees in order to be able to outperform the standard $\ell_1$.

Besides the above concrete descriptions about what was done in the paper, we would like to add that the results here are a nice complement to a large collection of the results already established for both, partial and hidden partial modifications of the standard $\ell_1$. Namely, in our prior work we settled the asymptotic regime of these problems through the analysis of both, the phase-transition and the large deviations phenomena. These phenomena occur in large dimensional settings. Studying finite dimensional settings is typically very hard and much less is known about the algorithms' behavior in those regimes. It is here that we for the first time managed to capture the finite dimensional behavior of the partial $\ell_1$.

Finally, fairly routine adjustments of the techniques introduced here can be done so that they can handle various modification/extensions of the problems considered here. Some of them we will present in a couple of forthcoming papers.

\begin{singlespace}
\bibliographystyle{plain}
\bibliography{l1hidparfinn1Refs}
\end{singlespace}

\end{document}